\documentclass[11pt,epsf]{article}
\usepackage{graphics,color,boxedminipage,bm,amsmath}
\usepackage{epsfig}
\usepackage{epsfig}
\usepackage{enumerate}
\usepackage{wrapfig}
\usepackage{graphics}
\usepackage{verbatim}
\usepackage{alltt}
\newtheorem{conj}{Conjecture}

\begin{document}
\title{\Large \bf  A conditional greedy   algorithm for edge-coloring}
\author
{\normalsize Mark K. Goldberg \\
Department of Computer Science,\\
Rensselaer Polytechnic Institute \\
Troy, NY, 12180.\\
{\tt goldberg@cs.rpi.edu}
}
\date{}
\maketitle
\begin{abstract}
We present a novel algorithm for edge-coloring of multigraphs.
The correctness of this algorithm for multigraphs with $\chi' > \Delta +1$ 
($\chi'$ is the chromatic edge number and $\Delta$ is the maximum vertex 
degree) would prove a long standing conjecture in edge-coloring of
multigraphs.
\end{abstract}

\section{Introduction.}

The {\bf chromatic index} $\chi'(G)$ of a 
multigraph $G(V,E)$ is the minimal number of colors that can be assigned 
to the edges of $G$ so that no two adjacent edges receive the same color.
Clearly, $$ \Delta (G)\leq \chi'(G),$$
where $\Delta(G)$ is the maximal vertex degree in $G$.
The famous result by Vizing (\cite{V-65}) establishes 
$\chi' =\chi'(G) \leq \Delta(G) + p(G)$, where $p(G)$ is the
maximal number of parallel edges in $G$. 
In particular, for graphs, $\Delta(G) \leq \chi'(G)\leq \Delta(G) +1$.
A polynomial algorithm to determine the exact value of the chromatic index
is unlikely to exist, as
the problem was proved by Holyer (\cite{Ho-81}) to be NP-hard even 
for cubic graphs.
It is suspected that for every multigraph with $\chi' > \Delta +1$,
its chromatic index is determined by the parameter $\omega(G)$ called the
multigraph {\bf density}:
$$
\omega(G)=\max_{H \subseteq G}\lceil\frac{e(H)}{\lfloor v(H)/2\rfloor}\rceil ,
$$
where $H$ is a sub-multigraph of $G$,  and $v(H)$ (resp. $e(H)$) denotes the 
number of vertices (resp. edges) in $H$.
It is easy to see that $\omega(G) \leq \chi' (G)$ for every multigraph $G$. 
A conjecture connecting $\chi', \omega$, and $\Delta$ was independently 
proposed by 
Goldberg (\cite{G-73}), Gupta(\cite{Gupta}), and Seymour (\cite{S-79}) more 
than 40 years ago (see \cite{TheBook}): if $\chi' > \Delta+1$, then 
$\chi' = \omega$.  
In (\cite{G-84}), it was additionally conjectured that if $\Delta > \omega$, 
then $\chi' = \Delta$. 
The value $\max\{\Delta, \omega\}$ is the polynomially computable
fractional chromatic index of $G$ (\cite{C-Y-Z-09, Sch-2003, S-79}).
See  \cite{C-Y-Z-09, G-84, HaxKier, K-96, N-90, TheBook, T-00} 
for some partial results towards the proof.
All these results are obtained by using algorithms that are based on the 
recoloring of  maximal two-colored chains and the operation of fan-recoloring 
discovered by Vizing and by Gupta.
In this paper, we present a different type of an algorithm which colors the 
edges in a pre-determined
sequence so that every intermediate partial edge-coloring is 
{\bf admissible.} 
Proving that the algorithm colors all edges would imply  the conjecture above.\\

{\bf Notations}.
Given set $S \subseteq V(G)$, $G[S]$ denotes the subgraph induced on $S$. 
For $S,T \subseteq V(G)$, $deg(S,T)$ denotes the number of edges $xy$ 
such that $x \in S$ and $y \in T$.  
The degree of $x$ is  $deg(x) = deg(\{x\}, V - x);$
if $x \in S$, the degree of $x$ in $G[S]$ is denoted $deg_S(x)$;
$\delta_S(x)=deg(x) - deg_S(x)$;
$\delta(S) = deg(S, V- S)= \sum_{x \in S} \delta_S(x).$
We say that edge $e= xy$ {\bf crosses} set $S \subseteq V$, 
if $|\{x, y\} \cap S| = 1$.

The multigraphs considered in this paper are connected and have 
more than one vertex; thus, the degree of every vertex is positive. 
The colors used for edge-colorings are  always
integers from the interval $[1,k]$, for some $k \geq 1$; 
such colorings are called $k$-edge-colorings. Obviously, every 
$k$-edge-coloring is also a $(k+1)$-edge-coloring.

A coloring for which all edges are colored is called {\bf complete}.
A {\bf partial} edge-coloring may have edges that are not colored; 
it is still required that any two adjacent colored edges 
have distinct colors. If $\phi$ is a partial edge-coloring
and edge $e$ is not colored by $\phi$, then we write $\phi(e) = \emptyset$.
Given two partial $k$-colorings 
$\phi_1$ and $\phi_2$, we say that $\phi_2$ is an
extension of $\phi_1$, denoted $\phi_1 \preceq \phi_2$,
if for every edge $xy$ which is $\phi_1$-colored, 
$\phi_2(xy) = \phi_1(xy)$.
If coloring $\phi_2$ is obtained from $\phi_1$ by coloring an edge $e$ 
using color $i$, we  write $\phi_2 = ext(\phi_1, e, i)$.

Let $\phi$ be a partial $k$-edge-coloring of $G$ and $S \subseteq V(G)$.
Then
$E_{\phi}^{(un)}(S)$ denotes the set of all edges in $E(S)$ 
that are not colored. 
For $i \in [1,k]$, an edge is called $i$-{\bf free}, if it is not colored and 
it is not adjacent to any edge colored $i$. 
A vertex $x$ is called $i$-{\bf free},
if it is incident to an $i$-free edge.
The set of all $i$-free vertices in $S$ is denoted $S_{\phi}^{(i)}$.
The cover value of $S$ by $\phi$ is defined as
$$
cov_{\phi}(S) = \sum_{i = 1}^k\lfloor \frac{|S_{\phi}^{(i)}|}{2}\rfloor.
$$ 
A partial coloring $\phi$ is said to  cover set $S \subseteq V(G)$, if 
$cov_{\phi}(S) \geq |E_{\phi}^{(un)}(S)|.$ A partial coloring is called
{\bf admissible}, if it covers every subset $S \subseteq V$.

Further in the paper, $k = \omega(G)$. 

\section{Edge-coloring algorithm}

The algorithm Conditional\_Greedy presented below colors the edges 
one-by-one using the smallest
available colors so that the resulting partial colorings are admissible. 
The edges are colored in the order developed by procedure {\it Reorder}, 
which first sorts out the vertices and then the edges.  

{\bf Procedure Reorder $G$;}\\ {\bf Input:} A connected graph $G$;\\ 
{\bf Output:} An ordered sequence $\{e_1,\ldots, e_m\}$ of the edges in
$G$.
\begin{itemize}
\item let $|V(G)| = n$;\\
select $x_1$ as a vertex of the maximal degree in $G$;
\item
for every $i \in [2,n]$, select $x_i$ as a 
vertex $z\in V-\{x_1, \ldots, x_{i-1}\}$ which maximizes\\ 
$deg(\{z\}, \{x_1, \ldots, x_{i-1}\})$; 
\item
the edges are sorted lexicographically:
for any two edges $e_a =x_i x_j$ and $e_b = x_k x_l$, where 
$i < j$ and $k < l$, $e_a$ precedes $e_b$ if either $i < k$, 
or $i=k$ and $j < l$  (the equality for parallel edges is broken arbitrary).
\end{itemize}
 
{\bf Procedure Conditional\_Greedy;}\\
{\bf Input:} A connected multigraph $G$; colors $1,2, \ldots, \omega$;
{\bf Output:} An edge-coloring $\phi$.
\begin{enumerate}
\item
{\it Reorder}($G$); 
let $e_1, e_2, \ldots, e_m$ be the output order of the edges; 

\item
For $i =1, \ldots, m$, set $\phi(e_i) = \emptyset$;
\item
for  $i = 1, \ldots, m$, \\
\hspace*{0.2in} \{find the smallest $c \in [1, \omega]$,
for which $ext(\phi, e_i, c)$ is admissible;\\
\hspace*{0.3in}if such a color does not exist,\\
\hspace*{0.6in} halt and output $\phi$;\\
\hspace*{0.3in}else $\phi = ext(\phi, e_i, c);$\}
\item output $\phi$.
\end{enumerate}

\begin{conj}
If $\chi^{\prime}(G) > \Delta(G) +1$, then Conditional\_Greedy 
outputs a complete $\omega$-edge-coloring.
\end{conj}
{\bf Acknowledgment.} The author is grateful to Prof. B. Toft for help with 
references.
\bibliographystyle{abbrv}

\bibliography{pap.bib}

\begin{thebibliography}{10}

\bibitem{C-Y-Z-09}
G.~Chen, X.~Yu, and W.~Zang.
\newblock Approximating the chromatic index of multigraphs.
\newblock {\em Journal of Combinatorial Optimization}, 21(Issue 2):219--246,
  February 2011.

\bibitem{G-73}
M.~K. Goldberg.
\newblock On multigraphs of almost maximal chromatic class.
\newblock {\em Discret. Analiz}, 23:3--7, 1973.
\newblock In Russian.

\bibitem{G-84}
M.~K. Goldberg.
\newblock Edge coloring of multigraphs: Recoloring technique.
\newblock {\em J. Graph Theory}, 8:123--137, 1984.

\bibitem{Gupta}
B.~P. Gupta.
\newblock On the chromatic index and the cover index of a multigraph.
\newblock {\em In: Dold, A. and Eckmann, B. editors, Theory and Applications of
  Graphs}, (Proceedings, Michigan, May 11-15, 1976, volume 642 of Lecture Notes
  in Mathematics):204--215.
\newblock Springer, Berlin.

\bibitem{HaxKier}
P.~E. Haxell and H.~A. Kierstead.
\newblock Edge coloring multigraphs without small dense subsets.
\newblock {\em Discrete Mathematics}, 338:2502--2506, 2015.

\bibitem{Ho-81}
I.~Holyer.
\newblock The {NP}-completeness of edge coloring.
\newblock {\em SIAM J. Comput.}, 10:718--720, 1981.

\bibitem{K-96}
L.~Kahn.
\newblock Asymptotics of the chromatic index for multigraphs.
\newblock {\em J. Combin. Theory Ser. B}, 68:195--235, 1996.

\bibitem{N-90}
T.~Nishizeki and K.~Kashiwagi.
\newblock On the 1.1 edge-coloring of multigraphs.
\newblock {\em SIAM J. Discrete Math.}, 3:391--410, 1990.

\bibitem{Sch-2003}
A.~Schrijver.
\newblock Combinatorial optimization: Polyhedra and efficiency.
\newblock {\em Prentice Hall, Upper Saddle River, NJ}, A, 2003.

\bibitem{S-79}
P.~D. Seymour.
\newblock Some unsolved problems on one-factorizations of graphs.
\newblock {\em Graph Theory and Related Topics, Academic Press}, 1979.
\newblock Bondy and Murty, eds.

\bibitem{TheBook}
M.~Steibitz, D.~Schide, B.~Toft, and L.~M. Favrholt.
\newblock Graph {E}dge {C}oloring.
\newblock {\em WILEY Series in Discrete Mathematics and Optimization.}, 2012.

\bibitem{T-00}
V.~A. Tashkinov.
\newblock On an algorithm for the edge coloring of multigraphs.
\newblock {\em Discretn. Anal. Issled. Oper.}, 1(7):72--85, 2000.
\newblock in Russian.

\bibitem{V-65}
V.~G. Vizing.
\newblock Critical graphs with a given chromatic class.
\newblock {\em Discret. Analiz}, 5:9--17, 1965.
\newblock in Russian.

\end{thebibliography}
\end{document}